\documentclass[12pt,a4paper,reqno]{amsart}
\usepackage{amssymb,color,currfile,pdfpages,mathtools,xcolor}

\edef\restoreparindent{\parindent=\the\parindent\relax}
\usepackage{parskip}
\restoreparindent

\usepackage[english]{babel}
\usepackage{verbatim,here}
\usepackage[T1]{fontenc}
\usepackage{epstopdf}
\usepackage{floatflt,graphicx}
\usepackage{color,xcolor}
\usepackage{enumerate}
\usepackage{animate}
\usepackage{a4wide}
\usepackage{amsmath, amssymb}
\usepackage{amsthm}
\usepackage{float}
\usepackage{pdfpages}
\usepackage{orcidlink}
\usepackage{hyperref}

\newcounter{minutes}
\divide\time by 60
\newcounter{hours}
\multiply\time by 60 \addtocounter{minutes}{-\time}

\newenvironment{pf}[1][]{%
	\vskip 3mm
	\noindent
	\ifthenelse{\equal{#1}{}}%
	{{\slshape Proof. }}%
	{{\slshape #1.} }%
}%
{\qed\bigskip}

\dedicatory{}
\commby{}
\theoremstyle{plain}
\newtheorem{thm}[equation]{Theorem}
\newtheorem{cor}[equation]{Corollary}
\newtheorem{lem}[equation]{Lemma}
\newtheorem{example}[equation]{Example}

\theoremstyle{definition}

\theoremstyle{remark}
\newtheorem{rem}[equation]{Remark}
\newtheorem{conj}[equation]{Conjecture}

\numberwithin{equation}{section}

\newcommand{\beq}{\begin{equation}}
	\newcommand{\eeq}{\end{equation}}
\newcommand{\ben}{\begin{enumerate}}
	\newcommand{\een}{\end{enumerate}}
\newcommand{\bequu}{\begin{eqnarray*}}
	\newcommand{\eequu}{\end{eqnarray*}}
\newcommand{\bequ}{\begin{eqnarray}}
	\newcommand{\eequ}{\end{eqnarray}}

\newcommand{\Bn}{ {\mathbb{B}^n} }

\begin{document}
	\thispagestyle{empty}
	\def\thefootnote{}
	\newcommand{\re}{{\rm Re\,}}
	\newcommand{\capacity}{{\mathop{\mathrm{cap}}}}
	\newcommand{\capacitydenL}{{\mathop{\mathrm{cap\ \underline{dens}}}}}
	\newcommand{\capacitydenU}{{\mathop{\mathrm{cap\ \overline{dens}}}}}
	
	\newcommand{\co}{{\overline{\operatorname{co}}}}
	\newcommand{\T}{{\mathcal T}}
	\newcommand{\U}{{\mathcal U}}
	\newcommand{\es}{{\mathcal S}}
	\newcommand{\LU}{{\mathcal{LU}}}
	\newcommand{\ZF}{{\mathcal{ZF}}}
	\newcommand{\IR}{{\mathbb R}}
	\newcommand{\IN}{{\mathbb N}}
	\newcommand{\IC}{{\mathbb C}}
	\newcommand{\IT}{{\mathbb T}}
	\newcommand{\ID}{{\mathbb D}}
	\newcommand{\IB}{{\mathbb B}}
	\newcommand{\K}{{\mathcal K}}
	\newcommand{\X}{{\mathcal X}}
	\newcommand{\PP}{{\mathcal P}}
	\newcommand{\uhp}{{\mathbb H}}
	\newcommand{\Z}{{\mathbb Z}}
	\newcommand{\N}{{\mathcal N}}
	\newcommand{\M}{{\mathcal M}}
	\newcommand{\SCC}{{\mathcal{SCC}}}
	\newcommand{\CC}{{\mathcal C}}
	\newcommand{\st}{{\mathcal{SS}}}
	\newcommand{\D}{{\mathbb D}}
	\newcommand{\sphere}{{\widehat{\mathbb C}}}
	\newcommand{\image}{{\operatorname{Im}\,}}
	\newcommand{\Aut}{{\operatorname{Aut}}}
	\newcommand{\real}{{\operatorname{Re}\,}}
	\newcommand{\kernel}{{\operatorname{Ker}}}
	\newcommand{\ord}{{\operatorname{ord}}}
	\newcommand{\id}{{\operatorname{id}}}
	\newcommand{\mob}{{\text{\rm M\"{o}b}}}
	\newcommand{\Int}{{\operatorname{Int}\,}}
	\newcommand{\Sign}{{\operatorname{Sign}}}
	\newcommand{\inv}{^{-1}}
	\newcommand{\area}{{\operatorname{Area}}}
	\newcommand{\eit}{{e^{i\theta}}}
	\newcommand{\ucv}{{\operatorname{UCV}}}

	\def\be{\begin{equation}}
		\def\ee{\end{equation}}
	\newcommand{\sep}{\itemsep -0.01in}
	\newcommand{\seps}{\itemsep -0.02in}
	\newcommand{\sepss}{\itemsep -0.03in}
	\newcommand{\bee}{\begin{enumerate}}
		\newcommand{\eee}{\end{enumerate}}
	\newcommand{\pays}{\!\!\!\!}
	\newcommand{\pay}{\!\!\!}
	\newcommand{\blem}{\begin{lem}}
		\newcommand{\elem}{\end{lem}}
	\newcommand{\bthm}{\begin{thm}}
		\newcommand{\ethm}{\end{thm}}
	\newcommand{\bcor}{\begin{cor}}
		\newcommand{\ecor}{\end{cor}}
	\newcommand{\beg}{\begin{example}}
		\newcommand{\eeg}{\end{example}}
	\newcommand{\begs}{\begin{examples}}
		\newcommand{\eegs}{\end{examples}}
	\newcommand{\bdefe}{\begin{defin}}
		\newcommand{\edefe}{\end{defin}}
	\newcommand{\bprob}{\begin{prob}}
		\newcommand{\eprob}{\end{prob}}
	\newcommand{\bei}{\begin{itemize}}
		\newcommand{\eei}{\end{itemize}}
	
	\newcommand{\bcon}{\begin{conj}}
		\newcommand{\econ}{\end{conj}}
	\newcommand{\bcons}{\begin{conjs}}
		\newcommand{\econs}{\end{conjs}}
	\newcommand{\bprop}{\begin{propo}}
		\newcommand{\eprop}{\end{propo}}
	\newcommand{\br}{\begin{rem}}
		\newcommand{\er}{\end{rem}}
	\newcommand{\brs}{\begin{rems}}
		\newcommand{\ers}{\end{rems}}
	\newcommand{\bo}{\begin{obser}}
		\newcommand{\eo}{\end{obser}}
	\newcommand{\bos}{\begin{obsers}}
		\newcommand{\eos}{\end{obsers}}
	\newcommand{\bpf}{\begin{pf}}
		\newcommand{\epf}{\end{pf}}
	\newcommand{\ba}{\begin{array}}
		\newcommand{\ea}{\end{array}}
	\newcommand{\llra}{\longleftrightarrow}
	\newcommand{\lra}{\longrightarrow}
	\newcommand{\lla}{\longleftarrow}
	\newcommand{\Llra}{\Longleftrightarrow}
	\newcommand{\Lra}{\Longrightarrow}
	\newcommand{\Lla}{\Longleftarrow}
	\newcommand{\Ra}{\Rightarrow}
	\newcommand{\La}{\Leftarrow}
	\newcommand{\ra}{\rightarrow}
	\newcommand{\la}{\leftarrow}
	\newcommand{\ds}{\displaystyle}
	\newcommand{\psubset}{\subsetneq}
	
	\def\cc{\setcounter{equation}{0}   
		\setcounter{figure}{0}\setcounter{table}{0}}
	
	\def\cc{\setcounter{equation}{0}   
		\setcounter{figure}{0}\setcounter{table}{0}}
	

	\title[H\"older continuity and composition operators]
	{H\"older continuity and composition operators in pluriharmonic Bloch spaces}

	\author{Jie Huang}
	\address{Jie Huang \vskip0.05cm  School of Mathematical Sciences, Huaqiao
		University, Quanzhou, Fujian 362021,  P. R. China.
       \vskip0.025cm Department of Mathematics with Computer Science, Guangdong Technion - Israel
		Institute of Technology, Shantou, Guangdong 515063, P. R. China.}
	\email{jie.huang@hqu.edu.cn;
    jie.huang@gtiit.edu.cn}
	
	\author{Suman Das}
	\address{Suman Das\vskip0.05cm Department of Mathematics, Shantou University, Shantou, Guangdong, P. R. China.
		\vskip0.025cm Department of Mathematics with Computer Science, Guangdong Technion - Israel
		Institute of Technology, Shantou, Guangdong 515063, P. R. China.}
	\email{suman@stu.edu.cn; suman.das@gtiit.edu.cn}

	\author{Antti Rasila}
	\address{Antti Rasila \vskip0.05cm Department of Mathematics with Computer Science, Guangdong Technion - Israel
		Institute of Technology, Shantou, Guangdong 515063, P. R. China. \vskip0.025cm Department of Mathematics, Technion - Israel
		Institute of Technology, Haifa 3200003, Israel.}
	\email{antti.rasila@gtiit.edu.cn; antti.rasila@iki.fi}

	
	\begin{abstract}
		We show a H\"older estimate of order $1/n$ for pluriharmonic Bloch mappings in the unit ball $\Bn \subset \IC^n$ with respect to the Bergman metric. We apply this to obtain a sufficient condition  for the composition operator on the pluriharmonic Bloch space to be bounded below. As a partial converse, we also give a necessary condition for the boundedness below of the composition operator on the Bloch space of holomorphic mappings in $\Bn$.
	\end{abstract}
	

	\maketitle
	

	\section{Introduction and preliminaries}
	Let $\IC^n$	denote the complex Euclidean space of $n\ge 1$ complex variables. For $z=(z_1, z_2, \ldots, z_n) \in \IC^n$ and $w=(w_1, w_2, \ldots, w_n)\in \IC^n$, the inner product is defined as $$\langle   z, w \rangle \coloneqq \sum_{i=1}^{n}z_i\overline{w}_i,$$ and the norm is given by $$|z|\coloneqq \langle  z, z \rangle^{\frac{1}{2}}.$$ For $x\in \IC^n$ and $r>0$, let us write $B^n(x, r)\coloneqq\{z\in\IC^n: |z-x|<r\}$ and $S^{n-1}(x, r)\coloneqq \partial B^n(x, r)$. Let $B^n(r)\coloneqq B^n(0,r)$,
	$\partial B^n(r)\coloneqq\partial B^n(0,r)$,
	$\Bn\coloneqq B^n(1)$, and denote by $\partial\Bn$ the boundary of the unit ball. In the case of one complex variable, $\ID\coloneqq\mathbb B^1$ is the usual unit disk.
	
	The class of all holomorphic functions from $\IB^n$ into $\IC^n$ is denoted by $\mathcal A (\Bn, \IC^n)$. We denote by $\mathrm{Aut}(\Bn)$ the automorphism group consisting of all biholomorphic self-mappings of the unit ball $\Bn$. 
	
	\subsection {Pluriharmonic mappings} A twice continuously differentiable function $f$ on $\Bn$ is said to be \textit{pluriharmonic} if its restriction to every complex line is harmonic. This is equivalent to the condition that, for any $z\in \Bn$, 
	$$\frac{\partial ^2}{\partial z_i \partial \overline{z}_j}f(z)=0,\ \ \  \forall i, j =1, 2, \ldots, n.$$
	Denote by $\mathcal H (\Bn, \IC^n)$ the class of all pluriharmonic mappings from $\IB^n$ into $\IC^n$. Every pluriharmonic mapping $f\in \mathcal H (\Bn, \IC^n)$ admits a representation $f=h+\overline{g},$ where $h, g\in \mathcal A (\Bn, \IC^n)$, and the representation is unique if $g(0)=0$.
	
	For a pluriharmonic mapping $f$ in $\Bn$, let
	\[
	J_f(z)\coloneqq
	\begin{pmatrix}
		D h(z) & \overline{D g(z)} \\
		D g(z) & \overline{D h(z)}
	\end{pmatrix}
	\]
	be the Jacobian matrix of $f$, where $Dh(z)$ and $Dg(z)$ denote the complex Jacobian matrices of $h$ and $g$ at $z$. Let
	\[
	\Omega_f \coloneqq \{ z\in\mathbb{B}^n : \det Dh(z)\neq 0\}.
	\]
	
	\subsection{Bergman metric} For $z\in \Bn$, let 
	$$B(z)\coloneqq\frac{(1-|z|^2)\mathit{I}+\mathit{A(z)}}{(1-|z|^2)^2}$$
	be the Bergman matrix, where $\mathit{I}$ is the $n\times n$ identity matrix, and 
	\[
	\mathit{A(z)}\coloneqq
	\begin{pmatrix}
		z_1\overline{z}_1 & \cdots & z_1\overline{z}_n \\
		\vdots & \ddots &  \vdots   \\
		z_n\overline{z}_1 & \cdots &  z_n\overline{z}_n
	\end{pmatrix}.
	\]
	For a smooth curve $\gamma : [0, 1]\rightarrow \Bn$, let 
	$$\ell_{\beta} (\gamma)\coloneqq\int_{0}^{1} \langle  B(\gamma(t))\gamma'(t), \gamma'(t) \rangle^{1/2}\, dt.$$
	For $z,\, w\in \Bn$, define $$\beta(z, w) \coloneqq \inf_\gamma \ell_{\beta}(\gamma),$$ where the infimum is taken over all piecewise smooth rectifiable curves $\gamma$ connecting $z$ and $w$. We call $\beta$ the \textit{Bergman metric}, and note that this metric can be represented as (see \cite[Page~25]{Zhuk}) $$\beta(z, w)=\frac{1}{2}\log{\frac{1+|\varphi_z(w)|}{1-|\varphi_z(w)|}}, \quad z, w\in \Bn,$$
	where $\varphi_z(w)$ is the involutive automorphism that interchanges the points $0$ and $z$. In particular, for any \(w \in \mathbb{B}^n\), the automorphism \(\varphi_w: \mathbb{B}^n \to \mathbb{B}^n\) that interchanges 0 and $w$ is given by
	\begin{eqnarray}\label{autBn}
		\varphi_w(z) = \frac{w - P_w z - \sqrt{1-|w|^2} Q_w z}{1 - \langle z, w \rangle}, \quad z \in \mathbb{B}^n,
	\end{eqnarray}
	where \(P_w\) is the orthogonal projection of \(\mathbb{C}^n\) onto the subspace spanned by $w$, and \(Q_w \coloneqq I - P_w\) is the projection onto the orthogonal complement. It has the following useful property:
	\begin{eqnarray}\label{equiv1}
		1-|\varphi_w(z)|^2=\frac{(1-|z|^2)(1-|w|^2)}{|1-\langle z, w \rangle|^2}, \quad z \in \mathbb{B}^n.
	\end{eqnarray}
	
	For any $z,w\in\mathbb{B}^n$, the pseudo-hyperbolic distance between $z$ and $w$ is given by
	$\rho(z,w):=|\varphi_z(w)|.$ It is well-known that the pseudo-hyperbolic distance is invariant under automorphisms of $\mathbb{B}^n$, i.e., $\rho(g(z),g(w))=\rho(z,w)$ for all $g\in\text{Aut}(\mathbb{B}^n)$. By the definition of the Bergman metric, we get, for all $n\geq1$,
	$$\tanh\beta(z,w)=\rho(z,w), \quad \forall z,w\in\mathbb{B}^n.$$
	
	In particular, if $n=1$ (i.e., if $z,w\in\mathbb{D}$), the automorphism has the explicit form
	\begin{equation}\label{def2}
		\varphi_w(z)=\frac{w-z}{1-\overline{w}z},
	\end{equation} and the pseudo-hyperbolic distance reduces to the classical form
	$$\rho(z, w)=\left|\frac{z-w}{1-\overline{z} w}\right|.$$	
	In addition, a straightforward computation gives the useful identity
	$$1-\rho(z, w)^2=\frac{(1-|z|^2)(1-|w|^2)}{|1-\bar{z}w|^2}=(1-|w|^2)|\varphi_z'(w)|.$$

	\subsection{Bloch-type spaces on pluriharmonic mappings}
	For $f\in \mathcal{H}(\Bn, \IC^n)$, define $$D_f^{n, \alpha}(z)\coloneqq (1-|z|^2)^{\frac{\alpha(n+1)}{2n}}\Lambda_f(z)^{\frac{1}{n}},$$ where $\alpha>0$ and $$\Lambda_f(z)\coloneqq|\det Dh(z)|+|\det Dg(z)|.$$ The \textit{prenorm} $\|f\|_{P(n, \alpha)}$ is given by
	$$\|f\|_{P(n, \alpha)}\coloneqq\sup\limits_{z\in \Bn}D_f^{n, \alpha}(z).$$
	
	We call a mapping $f\in \mathcal{H}(\Bn, \IC^n)$ a \textit{Bloch-type mapping}, and write $f\in \mathcal{B}_{P(n, \alpha)}$, if 
	$$\|f\|_{P(n, \alpha)}<\infty.$$ 
	The quantity \(\|\cdot\|_{P(n,\alpha)}\) is not a norm, since for \(n\ge 2\) 
	the non‑zero mapping \(f(z)=(z_1,0,\dots,0)\) satisfies 
	\(\|f\|_{P(n,\alpha)}=0\). 
	We therefore keep the term \textit{prenorm}, and 
	\(\mathcal{B}_{P(n,\alpha)}\) is simply the set of pluriharmonic mappings 
	with finite prenorm. No linear structure is required on 
	\(\mathcal{B}_{P(n,\alpha)}\), and in this paper we work exclusively 
	with the prenorm \(\|\cdot\|_{P(n,1)}\).

	If $f \in \mathcal A (\Bn, \IC^n)$, i.e., $f$ is holomorphic, then $$\Lambda_f(z)=|\det Df(z)|,$$ and in this case, we write $f\in \mathcal{B}_{A(n, \alpha)}$ if $\|f\|_{P(n, \alpha)}<\infty.$ For example, $\mathcal{B}_{A(1,1)}$ is the classical Bloch space of analytic functions, and $\mathcal{B}_{A(1,\alpha)}$ is the $\alpha$-Bloch space. It is clear that $\mathcal{B}_{P(1, 1)}$ coincides with the harmonic Bloch space, which was studied by Colonna \cite{Colonna} as a generalization
	of the analytic Bloch space. For convenience, we denote the prenorm on $\mathcal{B}_{A(n, \alpha)}$ as $\|f\|_{A(n, \alpha)} := \|f\|_{P(n, \alpha)}$ for any holomorphic mapping $f \in \mathcal{A}(\mathbb{B}^n, \mathbb{C}^n)$. In particular, $\|f\|_{A(n,1)}$ denotes the prenorm of $f$ in the holomorphic Bloch space $\mathcal{B}_{A(n,1)}$. We refer to \cite{Anderson,Bonk,Chenshaolin-11,Pommerenke-70,Pommerenke} for more information on harmonic Bloch spaces. Some connections of Bloch spaces with related types of spaces can be found in \cite{CLR20, CPR1, CPR2, Samy-RIM}. More recently, the authors in \cite{DPQ, HPQ} investigated extreme points and support
	points of harmonic Bloch and $\alpha$-Bloch mappings.
	
	\subsection{Composition operators}
	For a holomorphic self-map $\varphi: \Bn\to \Bn$, the composition 
	operator $C_{\varphi}$ is defined by
	\begin{equation}\label{def1}
		C_{\varphi} f(z)\coloneqq f\circ\varphi(z),
	\end{equation}
	where $z\in\Bn$ and $f\in \mathcal H (\Bn, \IC^n)$. In particular, if $n=1$ and $h$ is analytic, then 
	$$|(C_{\varphi_w} h)'(0)|=(1-|w|^2)|h'(w)|.$$ 
	
	We say that $C_{\varphi}$ is \textit{bounded below} on 
	$\mathcal{B}_{P(n, \alpha)}$ if there is a constant $k>0$ such that 
	$$\|C_{\varphi}f\|_{P(n, \alpha)}\geq k\|f\|_{P(n, \alpha)},$$ 
	for all $f\in \mathcal{B}_{P(n, \alpha)}$.

	In \cite{Proceeding-AMS}, Ghatage, Yan, and Zheng showed that $C'_{\varphi_w} h$ is a Lipschitz function with respect to the pseudo-hyperbolic metric. The result, later sharpened by Xiong \cite{Xiong}, can be stated as follows.
	
	\begin{thm}\label{Thm-1}$($\cite{Proceeding-AMS,Xiong}$)$
		Let $h\in \mathcal{B}_{A(1, 1)}$. Then the inequality
		$$\left|D_h^{1, 1}(z_2)-D_h^{1, 1}(z_1)\right|\leq \frac{3\sqrt{3}}{2}\rho(z_1, z_2)\|h\|_{P(1, 1)},$$
		holds for all $z_1, z_2\in\ID$.
	\end{thm}
	Chen and Kalaj \cite{CD} extended the problem to the several dimensional setting, for holomorphic functions in $\IC^n$. In \cite{CHZ}, Chen, Hamada, and Zhu considered the case of complex-valued harmonic functions and showed that, remarkably, the sharp constant remains as $3\sqrt{3}/2$, which improved the result previously obtained by Huang, Rasila, and Zhu in \cite{HRZ}.  Further, they used this result to study the composition operators $C_{\varphi_w}$ on the harmonic Bloch-type spaces. 
	
	This paper is organized as follows. In Section 2, we consider the problem for pluriharmonic mappings and establish a H\"older estimate of order $1/n$ with respect to the Bergman metric, which reduces to a Lipschitz estimate when $n=1$. Then in Section 3, we consider composition operators on the pluriharmonic Bloch space and give a sufficient condition for the operators to be bounded below.
	
	\section{H\"older continuity of pluriharmonic Bloch-type mappings}
Here we show that $D_f^{n, 1}(z)$ satisfies a H\"older estimate of order $1/n$ with respect to the Bergman metric. In the special case $n=1$, this reduces to Lipschitz continuity. For our purpose, we recall two useful results of Chen, Ponnusamy, and Wang \cite{Chen'slemma}.
	\begin{lem}\label{Chen'slemma}\cite{Chen'slemma}
		For $x\in [0, 1]$, let 
		$$\phi(x)=x(1-x^2)^{\alpha(n+1)/2}\sqrt{\alpha(n+1)+1}\left[\frac{\alpha(n+1)+1}{\alpha(n+1)}\right]^{\alpha(n+1)/2}$$
		and 
		$$a_0(\alpha)=\frac{1}{\sqrt{\alpha(n+1)+1}}.$$
		Then $\phi(x)$ is increasing in $[0, a_0(\alpha)]$, decreasing in $[a_0(\alpha), 1]$, and $\phi(a_0(\alpha))=1$. 
	\end{lem}
	
	\begin{thm}\label{Chen'sthm}\cite{Chen'slemma}
		Let $h\in  \mathcal A (\Bn, \IC^n)$ such that $\|h\|_{P(n, \alpha)}=1$ and $\det Dh(0)=\lambda\in (0, 1]$. Then, for all $|z|\leq \frac{a_0(\alpha)+m_{\alpha}(\lambda)}{1+a_0(\alpha)m_{\alpha}(\lambda)}$, we have 
		\begin{eqnarray}\label{Cheneq1}
			|\det Dh(z)|\geq \re (\det Dh(z))\geq \frac{\lambda(m_{\alpha}(\lambda)-|z|)}{m_{\alpha}(\lambda)(1-m_{\alpha}(\lambda)|z|)^{\alpha(n+1)+1}},
		\end{eqnarray}
		where $\phi$ and $a_0(\alpha)$ are as in Lemma \ref{Chen'slemma}, and  $m_{\alpha}(\lambda)$ is the unique real root of the equation $\phi(x)=\lambda$ in the interval $[0, a_0(\alpha)]$. Moreover, for all $|z|\leq \frac{a_0(\alpha)-m_{\alpha}(\lambda)}{1-a_0(\alpha)m_{\alpha}(\lambda)}$, we have 
		\begin{eqnarray}\label{Cheneq2}
			|\det Dh(z)|\leq  \frac{\lambda(m_{\alpha}(\lambda)+|z|)}{m_{\alpha}(\lambda)(1+m_{\alpha}(\lambda)|z|)^{\alpha(n+1)+1}}.
		\end{eqnarray}
		Both inequalities \eqref{Cheneq1} and \eqref{Cheneq2} are sharp. 
	\end{thm}
	
	Based on ideas from the proof of Theorem \ref{Chen'sthm}, we can obtain the same lower bound for $\Lambda_f$ when $f$ is a pluriharmonic mapping. We include an argument here for the convenience of the reader.
	
	\begin{lem}\label{lem2.5}
		Let $f=h+\overline{g}\in \mathcal{H}(\Bn, \IC^n)$ satisfy $\|f\|_{P(n, \alpha)}=1$ and $$\det Dh(0)+\det Dg(0)=\lambda\in (0, 1].$$ Then, for all $|z|\leq \frac{a_0(\alpha)+m_{\alpha}(\lambda)}{1+a_0(\alpha)m_{\alpha}(\lambda)}$, we have 
		\begin{eqnarray}
			\Lambda_f(z)\geq \re (\det Dh(z)+\det Dg(z))\geq \frac{\lambda(m_{\alpha}(\lambda)-|z|)}{m_{\alpha}(\lambda)(1-m_{\alpha}(\lambda)|z|)^{\alpha(n+1)+1}},
		\end{eqnarray}
		where $a_0(\alpha)$ and  $m_{\alpha}(\lambda)$ are as in Theorem  \ref{Chen'sthm}.
	\end{lem}

	\begin{proof}
		Fix an arbitrary $\zeta\in\partial\mathbb{B}^n$ and define the one-variable holomorphic function
		\[
		F(u)\coloneqq \det Dh(u\zeta)+\det Dg(u\zeta),\quad u\in\mathbb{D},
		\]
		with $F(0)=\det Dh(0)+\det Dg(0)=\lambda\in(0,1]$.
		Because $\|f\|_{P(n,\alpha)}=1$, we have
		\[
		\Lambda_f(u\zeta)=|\det Dh(u\zeta)|+|\det Dg(u\zeta)|
		\le \bigl(1-|u|^2\bigr)^{-\frac{\alpha(n+1)}{2}},
		\]
		and consequently
		\[
		|F(u)|\le \bigl(1-|u|^2\bigr)^{-\frac{\alpha(n+1)}{2}},\qquad u\in\mathbb{D}.
		\]
		
		Set $a=m_{\alpha}(\lambda)$, $a_0=a_0(\alpha)$, and we consider two cases.
		
		\noindent\textbf{Case 1:} Suppose $0<\lambda<1$. Then $a<a_0$, and therefore, $0\in\Omega$.
		Define
		\[
		T(u)=(1-au)^{\alpha(n+1)}F(u),\quad u\in\mathbb{D}.
		\]
		Let $\Omega=\mathbb{D}_h(a,\operatorname{arctanh} a_0)$ be the hyperbolic disk centered at $a$ with hyperbolic radius $\operatorname{arctanh} a_0$. Since $a<a_0$, we have $\rho(a,0)=a<a_0$, so $0\in\Omega$.
		A direct estimate on the boundary $\partial\Omega$ shows
		\[
		|T(u)|\le\kappa\coloneqq\frac{\lambda}{a\sqrt{\alpha(n+1)+1}},
		\]
		and the maximum principle yields $|T(z)|\le\kappa$ for all $z\in\Omega$.
		
		Consider the M\"obius transformation
		\[
		G_{\lambda}(u)=\frac{\lambda(a-u)}{a(1-au)},
		\]
		which maps $\Omega$ conformally onto $\mathbb{D}(0,\kappa)$ and satisfies $G_{\lambda}(0)=\lambda$.
		Since $T(\Omega)\subset\overline{\mathbb{D}(0,\kappa)}=G_{\lambda}(\Omega)$ and $G_{\lambda}$ is univalent,
		there exists a holomorphic map $\omega:\Omega\to\Omega$ with $\omega(0)=0$ such that $T=G_{\lambda}\circ\omega$ on $\Omega$.
		
		Since $\Omega$ is a hyperbolic disk centered at the positive real point $a$ 
		with radius $\operatorname{arctanh} a_0$, it is a Euclidean disk symmetric with 
		respect to the real axis and containing $0$.  
		Take the unique conformal map $\Phi:\Omega\to\mathbb{D}$ normalized by 
		$\Phi(0)=0$ and $\Phi'(0)>0$.  
		By the Schwarz reflection principle, $\Phi$ satisfies 
		$$\Phi(\overline{z})=\overline{\Phi(z)},$$ hence its Taylor coefficients are 
		real and it maps $\Omega\cap\mathbb{R}$ onto $(-1,1)$.  
		Moreover, $\Phi'(x)>0$ on $\Omega\cap\mathbb{R}$, so $\Phi$ is strictly 
		increasing there. 
		Let $\Psi=\Phi^{-1}:\mathbb{D}\to\Omega$. Define
		\[
		t(w)=\frac{T(\Psi(w))}{\kappa},\qquad g(w)=\frac{G_{\lambda}(\Psi(w))}{\kappa},\qquad w\in\mathbb{D}.
		\]
		Then $t,g:\mathbb{D}\to\mathbb{D}$ are holomorphic, $t(0)=g(0)=\gamma\coloneqq\lambda/\kappa$, and $g$ is surjective.
		From $T=G_{\lambda}\circ\omega$ we obtain $t=g\circ\omega_1$ with $\omega_1=\Phi\circ\omega\circ\Psi:\mathbb{D}\to\mathbb{D}$ and $\omega_1(0)=0$.
		
		Now fix $u\in\left(0,\frac{a_0+a}{1+a_0a}\right]$. Because $\frac{a_0+a}{1+a_0a}<\frac{1}{a}$, we have $u\in\Omega\cap\mathbb{R}^+$.
		By the increasing property of $\Phi$, we let $w_0=\Phi(u)\ge0$.
		Applying Schwarz's lemma to $\omega_1$ gives $|\omega_1(w_0)|\le w_0$.
		
		We now consider the function $g$. Since $g$ is a M\"obius transformation mapping $\mathbb{D}$ onto $\mathbb{D}$, 
		the image of the closed disk $\{z:|z|\le w_0\}$ under $g$ is a closed disk 
		contained in $\overline{\mathbb{D}}$. 
		Because $g$ has real coefficients, this image disk is symmetric with respect 
		to the real axis. 
		Moreover, $g$ is strictly decreasing on $[-1,1]$,  and thus, on the 
		interval $[-w_0,w_0]$ the minimum of $\operatorname{Re} g$ is attained at 
		the right endpoint $z=w_0$. 
		Therefore, for every $z$ with $|z|\le w_0$, we have 
		$\operatorname{Re} g(z) \ge g(w_0)$.
		Applying this with $z=\omega_1(w_0)$ yields 
		$\operatorname{Re} t(w_0) \ge g(w_0)$.
		
		Therefore,
		\[
		\operatorname{Re} T(u)=\kappa\operatorname{Re} t(w_0)\ge\kappa g(w_0)=G_{\lambda}(u).
		\]
		Finally, for any $z\in\mathbb{B}^n$ with $|z|\le\frac{a_0+a}{1+a_0a}$, write $z=u\zeta$ where $u=|z|$.
		Then
		\begin{align*}
			\operatorname{Re}\bigl(\det Dh(z)+\det Dg(z)\bigr)
			&=\frac{\operatorname{Re} T(u)}{(1-au)^{\alpha(n+1)}}
			\ge\frac{G_{\lambda}(u)}{(1-au)^{\alpha(n+1)}}
			=\frac{\lambda(a-|z|)}{a(1-a|z|)^{\alpha(n+1)+1}}.
		\end{align*}
		Because $\Lambda_f(z)\ge|\det Dh(z)+\det Dg(z)|\ge\operatorname{Re}(\det Dh(z)+\det Dg(z))$, the desired inequality follows. 
		
	\textbf{Case 2:} If $\lambda=1$, then $a=a_0$ and $u_2:=\frac{a_0+a}{1+a_0a}= \frac{2a_0}{1+a_0^2}$.
		For any $\zeta\in\partial\mathbb{B}^n$, we define
		$F(u)=\det Dh(u\zeta)+\det Dg(u\zeta)$ as above. Then $F$ is holomorphic in
		$\mathbb{D}$ with $F(0)=1$ and
		$$|F(u)|\le(1-|u|^2)^{-\frac{\alpha(n+1)}{2}}.$$
		
		For $k\ge 2$, set $c_k = 1 - \frac{1}{k}$ and $F_k(u) = c_k F(u)$, then each $F_k$ is holomorphic in $\mathbb{D}$, and satisfies
		$$|F_k(u)|\le(1-|u|^2)^{-\frac{\alpha(n+1)}{2}},$$ with $F_k(0)=c_k\in(0,1)$.
		Applying Case~1 to $F_k$ gives, for all
		$|u| \le r_k := \frac{a_0 + m_{\alpha}(c_k)}{1 + a_0\,m_{\alpha}(c_k)}$,
		\[
		\operatorname{Re} F_k(u) \ge
		\frac{c_k\bigl(m_{\alpha}(c_k)-|u|\bigr)}
		{m_{\alpha}(c_k)\bigl(1-m_{\alpha}(c_k)|u|\bigr)^{\alpha(n+1)+1}} .
		\]
		Since $F_k = c_k F$, we have
		\[
		\operatorname{Re} F(u) \ge
		\frac{m_{\alpha}(c_k)-|u|}
		{m_{\alpha}(c_k)\bigl(1-m_{\alpha}(c_k)|u|\bigr)^{\alpha(n+1)+1}} .
		\]
		
		Now fix any $z$ with $0 < |z| \le u_2$ and set $u = |z|$,
		$\zeta = z/|z|$.  As $k\to\infty$, we have
		$c_k\to 1$, $m_{\alpha}(c_k)\to a_0$, and $r_k\to u_2$. Hence
		$u \le r_k$ for all sufficiently large $k$.  Letting $k\to\infty$, we have
		\[
		\operatorname{Re} F(u) \ge
		\frac{a_0 - |u|}{a_0(1 - a_0|u|)^{\alpha(n+1)+1}} .
		\]
		
		Therefore,
		\[
		\Lambda_f(z) \ge |F(u)| \ge \operatorname{Re} F(u)
		\ge \frac{a_0 - u}{a_0(1 - a_0 u)^{\alpha(n+1)+1}}
		= \frac{a_0 - |z|}{a_0(1 - a_0|z|)^{\alpha(n+1)+1}} .
		\]
		This completes the proof.
	\end{proof}

	We are now ready to establish the first main result of this paper. The following is an analogue of Theorem \ref{Thm-1} for pluriharmonic Bloch-type mappings in $\Bn$.
	
	
	\begin{thm}\label{thm1}
		Let $f\in \mathcal{B}_{P(n,1)}$. Then, for $z_1, z_2\in\Bn$,
		$$\left|D_f^{n,1}(z_2)-D_f^{n, 1}(z_1)\right|\leq M(n)\|f\|_{P(n, 1)}\big[\tanh \beta(z_1, z_2)\big]^{\frac{1}{n}},$$
		where $$M(n)=(n+2)^{\frac{1}{2n}}\left(\frac{n+2}{n+1}\right)^{\frac{n+1}{2n}}.$$
	\end{thm}		
	\begin{proof}
		Without loss of generality, we may assume that $\|f\|_{P(n, 1)}=1$ and $D_f^{n,1}(z_2)\leq D_f^{n, 1}(z_1)$. Let $\varphi\in \mathrm{Aut}(\Bn)$ be an automorphism such that $\varphi(0)=z_1$ and $\varphi^{-1}(z_2)=w$. Let us write $$F(z)=f\circ \varphi(z)=H(z)+\overline{G(z)},$$ where $H=h\circ \varphi$  and $G=g\circ \varphi$. It follows from \cite[Corollary 1.22]{Zhuk} that
		$$\tanh \beta(z_1, z_2)=\tanh \beta\left(\varphi^{-1}(z_1), \varphi^{-1}(z_2)\right)=\tanh \beta(0, w)=|w|.$$
		Since 
		$$|\det D\varphi(z)|=\left(\frac{1-|\varphi(z)|^2}{1-|z|^2}\right)^{(n+1)/2},$$
		we obtain that $$\|F\|_{P(n, 1)}=\|f\|_{P(n, 1)}=1.$$ A direct computation yields
		\begin{align*}
			\Lambda_F(0)^{1/n}&=\Lambda_f(\varphi(0))^{1/n}|\det D\varphi(0)|^{1/n}\\ &=(1-|z_1|^2)^{(n+1)/(2n)}\Lambda_f(z_1)^{1/n}\\ &=D_f^{n, 1}(z_1),
		\end{align*}
		and \begin{align*}
			D_f^{n, 1}(z_2)&=(1-|w|^2)^{(n+1)/(2n)}|\det D\varphi(w)|^{1/n}\Lambda_f(\varphi(w))^{1/n}\\ &=(1-|w|^2)^{(n+1)/(2n)}\Lambda_F(w)^{1/n}.\end{align*}
		Therefore, we find that 
		\begin{equation}\label{thm1eq}
			|D_f^{n,1}(z_2)-D_f^{n, 1}(z_1)|=\Lambda_F(0)^{1/n}-(1-|w|^2)^{(n+1)/(2n)}\Lambda_F(w)^{1/n}.
		\end{equation}
		
		If $\Lambda_F(0)=0$, then the result is obvious. Hence, we consider the case $\Lambda_F(0)\neq 0$, and we can assume that $\Lambda_F(0)=\lambda$. There exist $\theta_1, \theta_2\in [0, 2\pi]$ such that $$e^{i\theta_1}\det DH(0)+e^{i\theta_2}\det DG(0)=\Lambda_F(0)=\lambda.$$ 
		Applying Lemma \ref{lem2.5} to the function $$F_{\theta_1, \theta_2}(z):=e^{i\theta_1/n}H(z)+\overline{e^{i\theta_2/n}G(z)}$$ and $\alpha=1$, we have, for $|z|\leq \frac{a_0(1)+m_{1}(\lambda)}{1+a_0(1)m_{1}(\lambda)}$, 
		$$\Lambda_F(z)=\Lambda_{F_{\theta_1,\theta_2}}(z)\geq  \re \left(e^{i\theta_1}\det DH(z)+e^{i\theta_2}\det DG(z)\right)\geq \frac{\lambda(m_{1}(\lambda)-|z|)}{m_{1}(\lambda)(1-m_1(\lambda)|z|)^{n+2}},$$
		where \begin{align*}
			\lambda&=\sqrt{n+2}\left(\frac{n+2}{n+1}\right)^{(n+1)/2}m_1(\lambda)\left(1-m_1(\lambda)^2\right)^{(n+1)/2}\\ &= M(n)^{n}m_1(\lambda)\left(1-m_1(\lambda)^2\right)^{(n+1)/2},\end{align*}
		with $$M(n)=(n+2)^{\frac{1}{2n}}\left(\frac{n+2}{n+1}\right)^{\frac{n+1}{2n}}.$$
		\textbf{Case 1:} Let $|w| \leq m_1(\lambda)$. Since $m_1(\lambda)\leq \frac{a_0(1)+m_{1}(\lambda)}{1+a_0(1)m_{1}(\lambda)}$, it follows that
		$$(1-|w|^2)^{(n+1)/(2n)}\geq (1-m_1(\lambda)|w|)^{(n+1)/(2n)}\geq (1-m_1(\lambda)|w|)^{(n+2)/n}.$$
		By the concavity of the function $t^{1/n}$, we have
		$$(m_1(\lambda)-|w|)^{1/n}\geq m_1(\lambda)^{1/n}-|w|^{1/n}.$$ These,
		together with \eqref{thm1eq}, lead to
		\begin{align*}
			|D_f^{n,1}(z_2)-D_f^{n, 1}(z_1)|&=\lambda^{1/n}-(1-|w|^2)^{(n+1)/(2n)}\Lambda_F^{1/n}(w)\\
			&\leq \lambda^{1/n}-(1-|w|^2)^{(n+1)/(2n)}\frac{\lambda^{1/n}(m_1(\lambda)-|w|)^{1/n}}{m_1(\lambda)^{1/n}(1-m_1(\lambda)|w|)^{(n+2)/n}}\\
			&=\frac{\lambda^{1/n}}{m_1(\lambda)^{1/n}}\left(m_1(\lambda)^{1/n}-\frac{(1-|w|^2)^{(n+1)/(2n)}(m_1(\lambda)-|w|)^{1/n}}{(1-m_1(\lambda)|w|)^{(n+2)/n}}\right)\\
			&\leq \frac{\lambda^{1/n}}{m_1(\lambda)^{1/n}}|w|^{1/n}\\ & =M(n)\left(1-m_1(\lambda)^2\right)^{(n+1)/(2n)}|w|^{1/n}\\
			&\leq M(n)|w|^{1/n},
		\end{align*}
		which is the desired inequality.
		
		\noindent \textbf{Case 2:} Let $|w| > m_1(\lambda)$. Then 
		\begin{align*}
			|D_f^{n,1}(z_2)-D_f^{n, 1}(z_1)|&\leq \lambda^{1/n}\\&=M(n)m_1(\lambda)^{1/n}\left(1-m_1(\lambda)^2\right)^{(n+1)/(2n)}\\&<M(n)|w|^{1/n}.
		\end{align*}
		This completes the proof of the theorem.
	\end{proof}
	\begin{rem}\label{rem1}
		If $\|f\|_{P(n, 1)}=1$, then for every $\epsilon>0$ there exists a point $w\in \Bn$ such that $D_f^{n, 1}(w)>1-\epsilon$. This follows directly from the definition of the supremum. 
	\end{rem}

	As an interesting consequence of Theorem \ref{thm1}, we obtain the following result for planar harmonic Bloch functions.
	\begin{cor}
		Let $f\in \mathcal{B}_{P(1, 1)}$. Then, for almost every $z\in \ID$,
		\begin{align*}
			(1-|z|^2)\left(\left|\frac{\partial}{\partial z}D_f^{1, 1}(z)\right|+\left|\frac{\partial}{\partial \overline{z}}D_f^{1, 1}(z)\right|\right)\leq \frac{3\sqrt{3}}{2}\|f\|_{P(1,1)}.
		\end{align*}
		The constant $3\sqrt{3}/2$ on the right-hand side is the best possible.
	\end{cor}
	\begin{proof}
		Since $D_f^{1,1}$ is Lipschitz continuous with respect to the 
		pseudo-hyperbolic metric by Theorem~\ref{thm1}, it is locally 
		Lipschitz in the Euclidean metric on $\mathbb{D}$, and therefore 
		differentiable almost everywhere in $\mathbb{D}$ by Rademacher's 
		theorem. 
		
		Fix a point $z\in \mathbb{D}$ such that $D_f^{1,1}$ is differentiable. 
	For sufficiently small $r>0$ and $\theta\in[0,2\pi)$, set $w = z + r e^{i\theta}$.	
		It follows from Theorem~\ref{thm1} with $n=1$ that 
		\[
		|D_f^{1,1}(z) - D_f^{1,1}(w)| \le \frac{3\sqrt3}{2}\|f\|_{P(1,1)}\rho(z,w).
		\]
		Since $r$ and $\theta$ are arbitrary, this implies \begin{equation}\label{eq_cor}
			\max\limits_{\theta\in [0, 2\pi]}\left(\lim\limits_{r\rightarrow 0^{+}}\frac{|D_f^{1,1}(z)-D_f^{1,1}(w)|}{\rho(z, w)}\right) \le \frac{3\sqrt{3}}{2}\|f\|_{P(1,1)}.\end{equation}
		Now, we rewrite $$\max\limits_{\theta\in [0, 2\pi]}\left(\lim\limits_{r\rightarrow 0^{+}}\frac{|D_f^{1,1}(z)-D_f^{1,1}(w)|}{\rho(z, w)}\right)
		=	\max\limits_{\theta\in [0, 2\pi]}\left(\lim\limits_{r\rightarrow 0^{+}}\frac{|D_f^{1,1}(z)-D_f^{1,1}(w)|}{r}\cdot \frac{r}{\rho(z, w)}\right)$$ and observe that $$\lim\limits_{r\rightarrow 0^{+}} \frac{r}{\rho(z, w)} = \lim\limits_{r\rightarrow 0^{+}} \left|1-|z|^2 - \bar{z}r\eit\right| = 1-|z|^2.$$ On the other hand, $$\lim\limits_{r\rightarrow 0^{+}}\frac{D_f^{1,1}(w)-D_f^{1,1}(z)}{r} = e^{i\theta}\frac{\partial}{\partial z}D_f^{1, 1}(z)+e^{-i\theta}\frac{\partial}{\partial \bar{z}}D_f^{1, 1}(z)$$ is the directional derivative of $D_f^{1, 1}(z)$. Therefore, we have $$\max\limits_{\theta\in [0, 2\pi]}\left(\lim\limits_{r\rightarrow 0^{+}}\frac{|D_f^{1,1}(z)-D_f^{1,1}(w)|}{r}\right)
		= \left|\frac{\partial}{\partial z}D_f^{1, 1}(z)\right|+\left|\frac{\partial}{\partial \overline{z}}D_f^{1, 1}(z)\right|.$$ These, combined with \eqref{eq_cor}, lead to the desired inequality. Since almost every point is a differentiability point, the inequality holds for almost every $z$.
		
To verify that the constant $3\sqrt{3}/2$ is sharp, consider the analytic function
	\[
	f(z)=\frac{3\sqrt{3}}{4}z^{2}\in\mathcal{B}_{P(1,1)}.
	\]
	A direct computation gives
	\[
	D_f^{1,1}(z)=(1-|z|^2)|f'(z)|=\frac{3\sqrt{3}}{2}|z|(1-|z|^2),
	\]
	and thus $\|f\|_{P(1,1)}=1$. Taking $z=r>0$, we have
	$D_f^{1,1}(r)=\frac{3\sqrt{3}}{2}r(1-r^2)$, which is differentiable at $r>0$ with
\[
\frac{\partial}{\partial z}D_f^{1,1}(r)=\frac{\partial}{\partial \bar{z}}D_f^{1,1}(r)=\frac{3\sqrt{3}}{4}(1-3r^2).
\]
Therefore,
	\[
	(1-r^2)\left(\left|\frac{\partial}{\partial z}D_f^{1, 1}(r)\right|+\left|\frac{\partial}{\partial \overline{z}}D_f^{1, 1}(r)\right|\right)=\frac{3\sqrt{3}}{2}(1-r^2)(1-3r^2).
	\]
	Letting $r\to 0^{+}$, the right-hand side approaches $3\sqrt{3}/2$. Hence the constant cannot be improved.
	\end{proof}

	\section{Composition operators on the pluriharmonic Bloch space}
	Let $\varphi\colon \Bn\to \Bn$ be a holomorphic mapping. Let us define
	\begin{equation}\label{eq_new}
		\tau_{\varphi, n}(z)\coloneqq \left(\frac{1-|z|^2}{1-|\varphi(z)|^2}\right)^{\frac{n+1}{
				2n}}|\det D\varphi (z)|^{\frac{1}{n}} \quad (z\in \Bn).
	\end{equation}
	From \cite[Lemma 6]{Chenh}, it is known that $\tau_{\varphi, n}(z)\leq 1$. For the case $n=1$, Ghatage, Yan, and Zheng \cite{Proceeding-AMS} obtained the following sufficient condition such that the composition operator $C_{\varphi} $ is bounded below on the classical Bloch space $\mathcal{B}_{A(1,1)}$. 
	
	
	\begin{thm}\cite[Theorem 2]{Proceeding-AMS}
		Let $\varphi$ be an analytic self-mapping of the unit disk $\ID$. Suppose for some constants $r\in (0, 1/4)$ and $\epsilon>0$, for each point $w\in \D$, there is a point $z_w\in \ID$ such that 
		$$\rho(\varphi(z_w), w)<r \quad \text{and} \quad \tau_{\varphi, 1}(z_w)> \epsilon.$$
		Then $C_{\varphi} $ is bounded below on $\mathcal{B}_{A(1,1)}$.
	\end{thm}
	Here we extend the above result to the pluriharmonic Bloch space $\mathcal{B}_{P(n,1)}$. The following is the second main theorem of this paper.
	
	\begin{thm}\label{thm2}
		Let $\varphi: \Bn \to \Bn$ be a holomorphic mapping. Suppose that there are constants $r\in(0, 1/M(n))$ and $\epsilon>0$ such that for each $w\in \Bn$, there is a point $z_w\in \Bn$ satisfying $$\tanh \beta(\varphi(z_w), w)<r^n \quad \text{and} \quad \tau_{\varphi, n}(z_w)>\epsilon,$$ where $M(n)$ is defined as in Theorem \ref{thm1}. Then there exists a constant $k(n, r, \epsilon)>0$, depending only on $r, \epsilon$ and $n$, such that $\|C_\varphi f\|_{P(n,1)}\geq k(n, r, \epsilon)\|f\|_{P(n, 1)}$
        for all $f\in \mathcal{B}_{P(n, 1)}$.
	\end{thm}
	\begin{proof}
		Without loss of generality, we may assume that $\|f\|_{P(n, 1)}=1$. By Remark \ref{rem1}, there is a point $w\in \Bn$ such that $$D_f^{n, 1}(w)>1-\delta,$$ where $$\delta=a\left(1-M(n)r\right)$$ and $a\in (0, 1)$ is arbitrary. By the assumption, there is a point $z_w\in \Bn$ such that $$\big[\tanh \beta(\varphi(z_w), w)\big]^{\frac{1}{n}}<r<\frac{1}{M(n)}.$$ Therefore, we may apply Theorem \ref{thm1} to find that
		\begin{align*}
			D_f^{n, 1}(\varphi(z_w))&\geq D_f^{n, 1}(w)-M(n)\big[\tanh \beta(\varphi(z_w), w)\big]^{\frac{1}{n}}\\
			&> 1-\delta-M(n)r\\
			&=(1-a)(1-M(n)r)>0.
		\end{align*}
		Then, as $\tau_{\varphi, n}(z_w)>\epsilon$, we have 
		$$\|C_\varphi f\|_{P(n,1)}\geq D_f^{n, 1}(\varphi(z_w))\tau_{\varphi, n}(z_w)> (1-a)(1-M(n)r)\epsilon.$$
		Since $a\in (0, 1)$ is arbitrary, it follows that 
		$$\|C_\varphi f\|_{P(n,1)}\geq  k(n, r, \epsilon),$$
		where $$k(n, r, \epsilon)\coloneqq \left(1-M(n)r\right)\epsilon.$$ This completes the proof of the theorem.
	\end{proof}
	
	The following result is, in some sense, a partial converse of Theorem \ref{thm2}. It gives a necessary condition for the composition operator to be bounded below on the holomorphic Bloch space $\mathcal{B}_{A(n,1)}$. Since $\mathcal{B}_{A(n,1)}\subset\mathcal{B}_{P(n,1)}$ and the prenorms coincide on this subspace, the following theorem also provides the same necessary condition whenever $C_{\varphi}$ is bounded below on the pluriharmonic Bloch space $\mathcal{B}_{P(n,1)}$.
	
	\begin{thm}
		Suppose $C_{\varphi}$ is bounded below on $\mathcal{B}_{A(n,1)}$. Then there exist constants $\epsilon>0$ and $r \in (0,1)$ such that for every $w\in \Bn$, one can find a point $z_w\in \Bn$ satisfying
		$$\tanh \beta(\varphi(z_w), w)\leq r \quad \text{and} \quad \tau_{\varphi, n}(z_w)\geq \epsilon.$$
	\end{thm}
	\begin{proof}
		Since $C_{\varphi}$ is bounded below  on $\mathcal{B}_{A(n,1)}$, there is a constant $k>0$ such that for all $h\in \mathcal{B}_{A(n, 1)}$, the inequality
		\begin{eqnarray}\label{eqThm3(1)}
			\|C_{\varphi}h\|_{A(n, 1)}\geq k\|h\|_{A(n, 1)}
		\end{eqnarray}
		holds. Fix $w \in \Bn$ and consider  $\varphi_w \in \operatorname{Aut}(\mathbb{B}^n)$ as defined in \eqref{autBn}. It follows from \eqref{equiv1} that
		\begin{eqnarray*}
			\|\varphi_w\|_{A(n,1)}&=& \sup_{z \in \mathbb{B}^n}(1-|z|^2)^{\frac{n+1}{2n}}|\det D\varphi_w(z)|^{\frac{1}{n}}\\ &=&(1-|w|^2)^{\frac{n+1}{2n}}\sup_{z \in \mathbb{B}^n}\bigg(\frac{1-|z|^2}{|1-\langle z, w \rangle|^2}\bigg)^{\frac{n+1}{2n}}\\ \nonumber
			&=&(1-|w|^2)^{\frac{n+1}{2n}}\sup_{z \in \mathbb{B}^n}\bigg(\frac{1-|\varphi_w(z)|^2}{|1-|w|^2}\bigg)^{\frac{n+1}{2n}}=1.
		\end{eqnarray*}
		From (\ref{eqThm3(1)}) we get
		\[
		\|C_\varphi \varphi_w\|_{A(n,1)} \geq k	\|\varphi_w\|_{A(n,1)}=k.
		\]
		For $z\in\Bn$, we observe that
		\begin{eqnarray*}
			\|C_\varphi \varphi_w\|_{A(n,1)} &=& \sup_{z \in \mathbb{B}^n} (1-|z|^2)^{\frac{n+1}{2n}} \left|\det D(\varphi_w \circ \varphi)(z)\right|^{\frac{1}{n}}\\ \nonumber
			&=&\sup_{z \in \mathbb{B}^n}(1-|z|^2)^{\frac{n+1}{2n}}|\det D\varphi(z)|^{\frac{1}{n}}|\det D\varphi_w(\varphi(z))|^{\frac{1}{n}}\\ \nonumber
			&=&\sup_{z \in \mathbb{B}^n}(1-|z|^2)^{\frac{n+1}{2n}}|\det D\varphi(z)|^{\frac{1}{n}}\bigg(\frac{1-|w|^2}{|1-\langle \varphi(z), w \rangle|^2}\bigg)^{\frac{n+1}{2n}}\\ \nonumber
			&=&\sup_{z \in \mathbb{B}^n}(1-|z|^2)^{\frac{n+1}{2n}}|\det D\varphi(z)|^{\frac{1}{n}}\left(\frac{1-|\varphi_w(\varphi(z))|^2}{1-|\varphi(z)|^2}\right)^{\frac{n+1}{2n}}\\ \nonumber
			&=&\sup_{z \in \mathbb{B}^n}	(1-|\varphi_w(\varphi(z))|^2)^{\frac{n+1}{2n}}\tau_{\varphi, n}(z)\geq k.
		\end{eqnarray*}
		Consequently, there exists a point $z_w\in \Bn$ such that \begin{align*}
			\left(1-|\varphi_w(\varphi(z_w))|^2\right)^{\frac{n+1}{2n}}\tau_{\varphi, n}(z_w)\geq \frac{k}{2}.\end{align*}
		It follows that $\tau_{\varphi, n}(z_w)\geq k/2$ and $\left(1-|\varphi_w(\varphi(z_w))|^2\right)^{(n+1)/(2n)}\geq k/2$. Therefore,
		$$\tanh \beta(\varphi(z_w), w) = |\varphi_w (\varphi(z_w))|\leq \left[1-\left(\frac{k}{2}\right)^{\frac{2n}{n+1}}\right]^{\frac{1}{2}}.$$
		Now, we may choose $\epsilon=k/2$ and
		\[
		r=\left[1-\left(\frac{k}{2}\right)^{\frac{2n}{n+1}}\right]^{\frac{1}{2}} < 1.
		\]
		This completes the proof.
	\end{proof}

	\subsection*{Compliance of ethical standards}
	
	\subsection*{Author contributions}
	This research was conducted in Guangdong Technion under supervision of A. Rasila with J. Huang as the lead author. All authors contributed substantially to the paper and approved the final submitted version.
	
	\subsection*{Conflict of interests}
	The authors declare that there is no conflict of interest
	regarding the publication of this paper.
	
	\subsection*{Data availability statement} Data sharing not applicable to this article as no
	datasets were generated or analyzed during the current study.

	\subsection*{Acknowledgments}
	The authors thank Prof. Shaolin Chen for his important feedback on an earlier version of this paper. The authors also thank the anonymous referees for their helpful remarks. 
	
	\subsection*{Funding}
    The research was partially supported by the Natural Science Foundation of Guangdong Province (Grant no. 2024A1515010467) and the Li~Ka~Shing Foundation STU-GTIIT Joint Research Grant (No. 2024LKSFG06).

\end{document}